\def\N{{\bf N}}

\def\CC{{\rm\kern.24em
   \vrule width.02em
       height1.4ex depth-.05ex
   \kern-.26em C}}
\def\QQ{{\rm\kern.24em
   \vrule width.02em
       height1.4ex depth-.05ex
   \kern-.26em Q}}
\def\PP{{\rm I\kern-.25em P}}         \def\RR{{\rm I\kern-.25em R}}
\def\DD{{\rm I\kern-.25em D}}         \def\EE{{\rm I\kern-.25em E}}
\def\FF{{\rm I\kern-.25em F}}         \def\NN{{\rm I\kern-.23em N}}
\def\RRp{{\rm I\kern-.25em R}_{+}}
\def\IND{{\rm 1\kern-.25em I}}

\documentclass[11pt,a4paper]{article}
\usepackage{graphicx}
\begin{document}
\font\gothic = eufm10 \font\Bbb = msbm10

\begin{center}
    \Large{\bf  FUNCTIONAL LIMIT THEOREMS FOR LEVY  PROCESSES AND THEIR ALMOST SURE VERSIONS}
\end{center}

\bigskip

\centerline{By ELENA PERMIAKOVA (Kazan)}

Chebotarev Inst. of Mathematics and Mechanics, Kazan State
University

Universitetskaya 17, 420008 Kazan

e-mail: epermiakova@mail.ru

\medskip

{\small {\bf Abstract.} In the paper it proved a criterion of convergences in distribution in Skorohod space.
This criterion is applied to some special Levy processes.
The almost sure version of  limit theorems for this processes are obtained.}



\medskip
\section{Introduction.}

Consider the
 sequence of the random processes

\begin{equation}    \label{1}
X_n(x)=\frac{V(s_nx)}{a_n}-xb_n, \ \ x\in [0,1], \ \ n \in \NN
\end{equation}
where $s_n \to \infty \mbox{ as }\  n \to \infty$, $V=V(x)$, $x \in [0,\infty)$ is a Levy process and  $a_n, \ b_n \in \RR$. In this paper we  solve the question of approximation  $\alpha$-stable processes by Levy processes of type $(\ref{1})$.

 In Section 2 we prove the criterion  of convergence of sequence Levy processes $V_n$ in
distribution  in the Skorohod space $D[0,1]$
(Theorem 1). As a corollary we obtain a theorem of the convergence
in distribution of compound Poisson processes. This corollary can be used from approximation the Levy process by random sums with the Poisson index summation.

In Section 3  we  
introduce  conditions for the process
$V$ under which the sequence $X_n$ converges in distribution in
$D[0,1]$ to  a $\alpha$-stable  Levy process $Y_\alpha$, $0<\alpha<2$ (Theorem 2) and  conditions for process $V$ under
which $X_n$ converges in distribution in $D[0,1]$ to a
Wiener random process $W$ (Theorem 4).
In particular we consider $a_n=n^{1/{\alpha}}$ and  specify  conditions for the process $V$ under which $X_n$ converges in distribution to $Y_\alpha$, $0<\alpha<2$ (Theorem 3).

In Section 4 it proved almost sure versions of the limit theorems
from Section 3. We describe the sequences $(s_n)$ such 
that the sequence of measures
\begin{equation}   \label{2}
Q_n(\omega)=\frac{1}{\ln n}\sum_{k=1}^{n}\frac{1}{k}\delta_{X_k(t)(\omega)}, \ \ n\in \NN
\end{equation}
where $\delta_x$ denotes the measure of unit mass, concentrated in
the point $x$, converges weakly to a distribution of $Y_\alpha$ or  $W$ in
$D[0,1]$ for almost all $\omega \in \Omega$.

 The
proofs of our results are based on the criterion for integral type
almost sure version of a limit theorem which was obtained in
Chuprunov and Fazekas [4].
\bigskip

\section{Functional limit theorems.}

\medskip
We will denote by $\stackrel{d}{\to}$ the convergence in
distribution, by $\stackrel{d}{=}$ the equality in distribution, by $\phi_{\xi}$ the characteristic function of the random variable $\xi$. We will suppose $\sum_{i\in\emptyset}a_i=0.$

  Recall, that $V(t)$, $t\ge 0$ is a Levy process if $V(0)=0$ and $V$ is
a homogeneous stochastically continuous random process with
independent increments, such that it's trajectories belong to the
Skorohod space. In the paper we will denote by the same symbols
the random process and the random element corresponding to this
random process. 

Using the Levy's representation (see [5], sect. 18) we can assume
that the characteristic function of the Levy  process $V=V(\gamma, \sigma, L,R)$ is

$$                                         \label{Levy}
\varphi_{V(t)}(x) =  {\EE}\left(e^{ixV(t)}\right) =
\psi\left(t, x, L(y), R(y), \sigma , \gamma
\right)=\exp{\left[itx \gamma -\right.}
$$

$$
\left.-\frac{\sigma^2x^2t}{2}+t\left\{
\int\limits_{-\infty}^0\left(e^{ixy}-1-\frac{ixy}{1+y^2}\right)dL(y)+
\int\limits^{\infty}_0\left(e^{ixy}-1-\frac{ixy}{1+y^2}\right)dR(y)\right\}\right]\,,
$$
$x\in{\RR}$

Here $L(y)$ is (left-continuous and) non-decreasing on $(-\infty, 0)$ with $L(-\infty)=0$ and
$R(y)$ is (right-continuous and) non-decreasing on $(0, \infty)$ with $R(\infty)=0$ and 
\begin{equation}               \label{usl}
\int_{-\varepsilon}^0y^2dL(y)+\int_0^{\varepsilon}y^2dR(y) <+\infty
\mbox{ for all } \varepsilon>0.
\end{equation}
 Not bounded the community, we can assume $\gamma=0$ and we will denote $V(\sigma,L,R)=V(0,\sigma,L,R) $.

The main result of the section is the following theorem.

{\bf Theorem 1. }{\it Let $V_n=V_n( \sigma_n, L_n,R_n), \ V=V( \sigma, L,R)$ be Levy processes. Then
\begin{equation}      \label{con2}
 V_n\stackrel{d}{\to } V \mbox{ as}\  n \to \infty
\end{equation}
in $D[0,1]$
if and only if
\begin{equation}   \label{con1}
V_n(1)\stackrel{d}{\to}V(1) \  \mbox{ as }\  n \to \infty.
\end{equation}
 }

{\bf Proof. }
The proof of the necessity is evident. We will prove the sufficiency.
By Theorem 2, sect.19 in [5],  ($\ref{con1}$) implies 
$$ L_n(y)\to L(y), \ \ R_n(y)\to R(y),\ \ {\mbox as } \ \ n \to \infty$$
and
$$\lim_{\varepsilon \to 0}\overline{\lim_{n\to\infty}}\left\{ \int_{-\varepsilon}^{0}u^2dL_n(u)+\sigma_n^2+\int_{0}^{\varepsilon}u^2dR_n(u) \right\}=$$
$$\lim_{\varepsilon \to 0}\underline{\lim_{n\to\infty}}\left\{ \int_{-\varepsilon}^{0}u^2dL_n(u)+\sigma_n^2+\int_{0}^{\varepsilon}u^2dR_n(u) \right\}=\sigma^2.$$
Therefore by  Theorem 2 of sect. 19 in [5]  $V_n(t)\ \stackrel{d}{\to}V(t)$ as $n\to \infty$ for all $t\in [0,1].$ 

Let $0\le t_0<t_1<...<t_n\le 1.$
Introduce the notation $\Delta V_{ni}=V_n(t_i)-V_n(t_{i-1})$ and $\Delta V_i=V(t_i)-V(t_{i-1}).$
Since $\Delta V_{ni}, \ 1\le i\le k$ are independent random variables and $\Delta V_{ni}\stackrel{d}{\to}\Delta V_i$ as $n \to \infty$, we have:
$$(\Delta V_{n1},...,\Delta V_{nk})\stackrel{d}{\to}(\Delta V_1,...,\Delta V_k) \mbox{ as } n \to \infty. $$
Consequently, the finite dimensional distributions of $V_n$ converge to the finite dimensional distributions of $V$.

By Theorem 1 in [11] (p. 229), it's sufficiently to demonstrate that for all $\varepsilon>0$
\begin{equation}    \label{relcomp}
\lim_{h\to\infty}\overline{\lim_{n\to\infty}}\sup_{|t'-t''|\le h}P\{|V_n(t')-V_n(t'')|>\varepsilon\}. 
\end{equation}

By Theorem 4 in [12] the set of measures in $D[0,1]$,  corresponding to the processes  $V_n$ is relatively compact. Then by Theorem 1 in [12] we obtain $(\ref{relcomp})$.
Proof is completed.

 Let $\xi_{ni}, \xi_i \ i \in \N, $ be  independent identically distributed random variables for all $n \in\NN$, $\pi(t)$ is a Poisson random process with the intensity 1, $\xi_{ni}$ and $\pi(t)$ are independent. We will consider the Levy processes 
 \begin{equation}  \label{defVnprim}
 V_n'(t)=\sum_{i=1}^{\pi(t)}\xi_{ni}, \ t\in [0,1] 
 \end{equation}
 and
\begin{equation}    \label{defVprim}
V'(t)=\sum_{i=1}^{\pi(t)}\xi_{i}, \ \ t\in [0,1].
\end{equation}

   Such processes are called the compound processes.

{\bf Corollary 1. }{\it Let $V_n'$ and $V'$ be defined by $(\ref{defVnprim})$, $(\ref{defVprim})$. Suppose that
 \begin{equation}      \label{corollary}
 \xi_{ni}\stackrel{d}{\to}\xi_{i} \mbox{ as } n \to \infty.
 \end{equation}
  Then we have
  \begin{equation}             \label{col12}
  V_n' \stackrel{d}{\to}V' \mbox{ as } n\to \infty
   \end{equation}
   in $D[0,1]$. }

{\bf Proof.}
Observe that  ($\ref{corollary}$) implies
 $$\phi_{\xi_{n1}}(x)\to\phi_{\xi_1}(x) \mbox{ as } n \to \infty \ \ \mbox{ for all } x \in \RR. $$
 Since
 $$\phi_{V_n'(1)}(x)=\exp{\{\phi_{\xi_{n1}}(x)-1\}},\ \  x \in \RR$$
and
 $$ \phi_{V'(1)}(x)=\exp{\{\phi_{\xi}(x)-1\}}, \ \ x \in \RR, $$
 we have $$ \phi_{V_n'(1)}(x)\to\phi_{V'(1)}(x) \mbox{ as } n \to \infty \ \mbox{\rm for all } x \in \RR.$$
 Consequantly, $$ V_n'(1) \stackrel{d}{\to}V'(1) \mbox{ as } n\to \infty$$ and we can apply Theorem 1.
 Proof is completed.
 
 Consider now the random processes 
 \begin{equation}   \label{defVn2prim}
 V_n''(t)=\sum_{i=1}^{\pi(k_nt)}\xi_{ni}, \ \ t\in [0,1], \ \ k_n \in \N, \ \ n\in \NN.
 \end{equation}

 {\bf Corollary 2. }{\it Let $V_n''$ be defined by $(\ref{defVn2prim})$. Suppose that
 \begin{equation}      \label{corollary2}
 \sum_{i=1}^{k_n}\xi_{ni}\stackrel{d}{\to}\gamma \mbox{ as } n \to \infty,
 \end{equation}
 where $\gamma$ is an infinitely divisible random variable.
  Then one has
  \begin{equation}             \label{corollary12}
  V_n'' \stackrel{d}{\to}V'' \mbox{ as } n\to \infty
   \end{equation}
   in $D[0,1]$, where  $V''(t), \ t>0$ is a Levy process such that $V''(1)\stackrel{d}{=}\gamma$. }

{\bf Proof. }
The convergence ($\ref{corollary2}$) implies, that for all $x \in \RR$ $$\phi_{\sum_{i=1}^{k_n}\xi_{ni}}(x)\to\phi_{\gamma}(x) \mbox{ as } n \to \infty. $$
 The characteristic function of $V_n''(1)$ is
 $$\phi_{V_n''(1)}(x)=\exp{\{k_n(\phi_{\xi_{n1}}(x)-1)\}}, \ \ x\in \RR.$$
 Observe, that 
 $$k_n\ln\phi_{\xi_{ni}}(x)=k_n\ln(1-(1-\phi_{\xi_{ni}}(x)))=(\phi_{\xi_{ni}}(x)-1)k_n+k_no(\phi_{\xi_{ni}}(x)-1). $$
 Therefore $$ \phi_{V_n''(1)}(x)\to\phi_{V''(1)}(x) \mbox{ as } n \to \infty \mbox{ for all } x \in \RR.$$
 Consequantly, $$ V_n''(1) \stackrel{d}{\to}V''(1) \mbox{ as } n\to \infty$$ and we can apply Theorem 1.
 Proof is completed.

\section{The convergence to $\alpha$-stable processes.}
Let $0<\alpha<2$.

Consider  the conditions:
\begin{equation}                                         \label{4}
{L(-y)}/{|R(y)|} \to  {c_1}/{c_2}, \,\,\,\, \mbox{ as } \,\,\, y\to\infty\,,
\end{equation}

\begin{equation}                                         \label{5}
\frac{L(-y)+|R(y)|}{L(-yx)+|R(yx)|}  \to   x^\alpha, \,\,\,\,  \mbox{ as }  \,\,\, y\to\infty\,,
\end{equation}
where $c_1, c_2 \ge 0, \ c_1+c_2>0$ are valid.

By Theorem 2, sect.35 in [5], the conditions $(\ref{4}), (\ref{5})$ imply that an infinitely divisible random variable with functions $L(y)$ and $R(y)$ in Levy's representation of it's characteristic function belongs to the domain of attraction of the $\alpha$-stable law  having Levy's representation $L_\alpha(y)=\frac{c_1}{|y|^{\alpha}}, \ R_\alpha(y)=\frac{c_2}{-y^{\alpha}}$. Than it exists $a_n, \ s_n$ such that $s_n, \ a_n \to \infty$ as $n\to\infty$ and
\begin{equation} \label{conditions}
\mbox{ for all } y < 0 \ \ s_nL(a_ny)\to \frac{c_1}{|y|^{\alpha}}, \mbox{ as } n \to \infty,\ \ \mbox{ for all } y > 0\ \  s_n R(a_ny)\to \frac{c_2}{-y^{\alpha}} \mbox{ as } n \to \infty.
\end{equation}
 Consider the sequence of the random processes 
 \begin{equation}   \label{Xnmen2}
 X_n(t)=\frac{V(s_nt)}{a_n}-tb_n, \ \ t\in [0,1], \ \ n \in \NN,
 \end{equation}    
where $s_n, \ a_n\in \RR $ are defined by conditions $(\ref{conditions})$ and 
\begin{equation}    \label{bn}
b_n=\int_{-\infty}^{0} \frac{ z^3(1-a_n^2)}{(1+z^2)(1+a_n^2z^2)}ds_nL(a_nz)+ \int_{0}^{+\infty} \frac{ z^3(1-a_n^2)}{(1+z^2)(1+a_n^2z^2)}ds_nR(a_nz).
\end{equation}

Denote $Y_\alpha=V(0,\frac{c_1}{|y|^{\alpha}},\frac{c_2}{-y^{\alpha}})$. Then $Y_\alpha$ is a $\alpha$-stable process $(0<\alpha< 2)$.
Now we can formulate the next theorem of convergence to a $\alpha$-stable Levy processes.

\medskip
{\bf Theorem 2. } {\it Let ($\ref{4}$)  and ($\ref{5}$) be valid, $s_n$ and $a_n$ be satisfy $(\ref{conditions}$), $X_n$ be defined by ($\ref{Xnmen2}$), $b_n$ be defined in $(\ref{bn})$ and $\sigma=0.$ Then we have
$$
X_n\stackrel{d}{\to} Y_\alpha, \,\,\,{\it as}\,\,\,n\to\infty
$$
in $D[0,1]$.}

{\bf Proof. }
 The characteristic function of $X_n$ is
 $$\varphi_{X_n(t)}(x)=$$
 $$=\exp\left[t\left\{\int\limits_{-\infty}^{0}\left(e^{ixy}-1-\frac{ixy}{1+y^2}\right)d(s_nL(a_ny))+\int\limits_{0}^{\infty}\left(e^{ixy}-1-\frac{ixy}{1+y^2}\right)d(s_nR(a_ny))\right\}\right].$$
 Then $$X_n(1)\stackrel{d}{\to}Y_\alpha(1) $$ as $n \to \infty$,
and by Theorem 1 we obtain the affirmation of theorem.
 Proof is completed.

 Consider the functions
 \begin{equation}   \label{Theorem4L}
 L(x)=(c_1+\alpha_1(x))\frac{1}{|x|^{\alpha}}, \ \ \ x<0,
 \end{equation}
 \begin{equation}   \label{Theorem4R}
 R(x)=(c_2+\alpha_2(x))\frac{1}{-x^{\alpha}}, \ \ \ x>0,
 \end{equation}
 where  $\alpha_1(x)$ and $\alpha_2(x)$ are the functions such that 
 \begin{equation}    \label{alpha}
 \lim_{x\to-\infty}\alpha_1(x)=\lim_{x\to+\infty}\alpha_2(x)=0. 
 \end{equation}

 By Theorem 5, sect. 35 in [5], the infinity divisible random variable $\gamma$ belongs to the domain of normal attraction of the $\alpha$-stable low if and only if the functions $L(x)$ and $R(x)$  in  Levy representation of its characteristic function are $(\ref{Theorem4L})$ and $(\ref{Theorem4R})$, relatively.
 
  Then  
 \begin{equation}  \label{normattr}
 \mbox{\rm for all } x <0 \ \ nL(n^{1/{\alpha}}x)\to \frac{c_1}{|x|^{\alpha}},\ \  \mbox{\rm for all } x >0 \ \ nR(n^{1/{\alpha}}x)\to \frac{c_2}{-x^{\alpha}}\ \  \mbox{\rm as } n \to \infty.
 \end{equation}
 
 Consider the sequence of the random processes
 \begin{equation}   \label{xnprim}
 X_n'(x)=\frac{V(nt)}{n^{1/{\alpha}}} -tb_n',\ \ x\in [0,1], n \in \NN,
 \end{equation}
 where 
 \begin{equation}    \label{bnprim}
b_n'=\int_{-\infty}^{0} \frac{ z^3(1-n^{2/\alpha})}{(1+z^2)(1+n^{2/\alpha}z^2)}dnL(n^{1/{\alpha}}z)+ \int_{0}^{+\infty} \frac{ z^3(1-n^{2/\alpha})}{(1+z^2)(1+n^{2/\alpha}z^2)}dnR(n^{1/{\alpha}}z).
\end{equation}

Then the next theorem is valid

{\bf Theorem 3. }{\it Let $X_n'$ be defined by $(\ref{xnprim})$, $b_n$ be defined by $(\ref{bnprim})$ and $\sigma=0$. Then the conditions $(\ref{Theorem4L})$, $(\ref{Theorem4R})$, $(\ref{alpha})$ imply the convergence
\begin{equation}  \label{Theorem3}
X_n'\stackrel{d}{\to}Y_\alpha, \mbox{ as } n \to \infty 
\end{equation}
in $D[0,1]$.}

{\bf Proof. }
 The characteristic function of $X_n'$ is
 $$\varphi_{X_n'(t)}(x)=$$
 $$=\exp\left[t\left\{\int\limits_{-\infty}^{0}\left(e^{ixy}-1-\frac{ixy}{1+y^2}\right)dnL(n^{1/{\alpha}}y)+\int\limits_{0}^{\infty}\left(e^{ixy}-1-\frac{ixy}{1+y^2}\right)dnR(n^{1/{\alpha}}y)\right\}\right].$$
 Then $$X_n'(1)\stackrel{d}{\to}Y_\alpha(1) $$ as $n \to \infty$,
and by Theorem 1 we obtain $(\ref{Theorem3})$.
 Proof is completed.
 
 Let now $\alpha=2$.
 Consider the condition
  \begin{equation}                                         \label{9}
\frac{x^2\left(\int_{-\infty}^{-x}dL(y)+\int_{x}^{+\infty}dR(y) \right)}{\int_{-x}^{0}y^2dL(y)+\int_{0}^{x}y^2dR(y)}\to 0 \mbox{ as } x\to \infty.
\end{equation}
By Theorem 1, sect. 35 in [5], the condition ($\ref{9}$) is valid if and only if an infinitely divisible random variable $\gamma$ with functions $L(y)$ and $R(y)$ in Levy's representation of it characteristic function belongs to the domain of attraction of Gaussian law  having Levy's representation $\psi(t,x,0,0,\sigma_W).$
  Then exists $a_n$, $s_n$ such that $s_n, \ a_n \to \infty$ as $n\to\infty$ and
 \begin{equation}   \label{cond1Theor3}
 \mbox{ for all } y < 0\ \  s_nL(a_ny)\to 0, \ \ \mbox{ for all } y > 0\ \  s_nR(a_ny)\to 0, 
 \end{equation}
 $$\lim_{\varepsilon \to 0}\overline{\lim_{n\to\infty}}\left\{ \int_{-\varepsilon}^{0}u^2ds_n L(a_n u)+\frac{\sigma^2}{a_n^2}s_n+\int_{0}^{\varepsilon}u^2ds_n R(a_n u) \right\}=$$
  \begin{equation}   \label{cond2Theor3}
  \lim_{\varepsilon \to 0}\underline{\lim_{n\to\infty}}\left\{ \int_{-\varepsilon}^{0}u^2ds_n L(a_nu)+\frac{\sigma^2}{a_n^2}s_n+\int_{0}^{\varepsilon}u^2ds_n R(a_nu) \right\}=\sigma_W^2.
  \end{equation}
 
  Consider the sequence of the random processes 
 \begin{equation}   \label{Xnravno2}
 X_n(t)=\frac{V(s_nt)}{a_n}-tb_n, \ \ t\in [0,1], \ \ n \in \NN,
 \end{equation}    
where $s_n, \ a_n\in \RR $ are defined by conditions $(\ref{cond1Theor3})$,  $(\ref{cond2Theor3})$ and 
$$b_n=\int_{-\infty}^{0} \frac{ z^3(1-a_n^2)}{(1+z^2)(1+a_n^2z^2)}ds_nL(a_nz)+ \int_{0}^{+\infty} \frac{ z^3(1-a_n^2)}{(1+z^2)(1+a_n^2z^2)}ds_nR(a_nz)- $$
\begin{equation}    \label{bnravno2}
-i\lim_{\varepsilon\to 0}\left\{\int\limits_{-\varepsilon}^0z^2ds_nL(a_nz)+\int\limits^{\varepsilon}_0 z^2ds_nR(a_nz)\right\}.
\end{equation}

 {\bf Theorem 4. } {\it Let  $a_n$, $s_n$ be satisfied ($\ref{cond1Theor3}$),($\ref{cond2Theor3}$) and $X_n$ be defined by ($\ref{Xnravno2}$), $b_n$ be defined by  ($\ref{bnravno2}$).
 Then the condition $(\ref{9})$
implies 
$$
X_n\stackrel{d}{\to}\sigma_W W, \,\,\,{\it as}\,\,\,n\to\infty,
$$
in $D[0,1]$ . }

{\bf Proof. }
 The characteristic function of $X_n(t)$ is:
 $$\varphi_{X_n(t)}(x)=\exp\left[ -s_nt\frac{\sigma^2}{2a_n^2}x^2 \right.+ $$
 $$\left.+t\left\{\int\limits_{-\infty}^{0}(e^{ixy}-1-\frac{ixy}{1+y^2})ds_nL(a_ny)+\int\limits_{0}^{\infty}(e^{ixy}-1-\frac{ixy}{1+y^2})ds_nR(a_ny)\right\}\right].$$
Then $$X_n(1)\stackrel{d}{\to}\sigma_WW(1) $$ as $n \to \infty$.
 Then by Theorem 1 $$ X_n\stackrel{d}{\to}\sigma_WW.$$
 Proof is completed.

Consider the random processes
 \begin{equation}            \label{xcon}
  X_t(x)=\frac{V(tx)}{f(t)}-xg(t), \ x\in [0,1], t\in(0, \infty),
 \end{equation}
 
  \begin{equation}          \label{gt}
  g(t)=\int\limits_{-\infty}^{0}\frac{v^3(f(t)^2-1)}{f(t)(1+v^2)(f(t)^2+v^2)}dtL(f(t)v)+\int\limits_{0}^{\infty}\frac{v^3(f(t)^2-1)}{f(t)(1+v^2)(f(t)^2+v^2)}dtR(f(t)v)
  \end{equation}

 \section{Almost sure versions.}

 Now we will consider the sequence of measures $Q_n(\omega)$, defined in ($\ref{2}$) and connected with the random processes $X_n$, defined in ($\ref{1}$) with the condition for sequence $s_n$:

 \bigskip
 
 (A) for some $\beta>0$ the sequence $\left(\frac{s_n}{n^{\beta}}\right)$ is  increasing as $n\to \infty$.

 \bigskip
 Let $Y_\alpha$ is a $\alpha$-stable  random process $(0<\alpha\le 2)$.
We will denote by $\stackrel{w}{\to}$ the weak convergence of measures, 
by $\mu_\alpha$   the distribution of the random element $Y_\alpha$ and by $\rho(x,y)$ the usual metric in $D[0,1]$, with which $D[0,1]$ is a separable complete space.

Let $\xi_{ni}$ are independent identically distributed random variables.
Consider the sums $S_n=\sum_{i=1}^{k_n}\xi_{ni}.$ Below we will used the next lemma:

{\bf Lemma 1. }{\it Let $S$ be an infinitely divisible random variable such that 
\begin{equation}   \label{sn}
S_n\stackrel{d}{\to}S \mbox{ as } n \to\infty.
\end{equation}
Then for all $C> 0$ the next conditions are valid:}
$$C_1=\sup_{n \in \N}k_nP\{|\xi_{n1}|\ge C\}<+\infty,$$
 $$C_2=\sup_{n \in \N}k_nE(\xi_{n1}I_{\{|\xi_1|<C\}})<+\infty, $$
$$C_3=\sup_{n\in \N}k_nD(\xi_{n1}I_{\{|\xi_{n1}|<C\}})<+\infty,$$ 
 
{\bf Proof.}
By Theorem 1, sect. 23 in [5], the convergence $(\ref{sn})$ implies that
$$\sup_{n\in\N}k_nE\frac{\xi_{n1}^2}{1+\xi_{n1}^2}<+\infty. $$
Observe, that
$$\sup_{n\in\N}k_nP\{|\xi_{n1}|\ge C\}=\sup_{n\in\N}k_nEI_{\{|\xi_{n1}|\ge C\}}=\sup_{n\in\N}k_n\frac{C^2+1}{C^2}E\left(\frac{C^2}{1+C^2}I_{\{|\xi_{n1}|\ge C\}}\right)\le $$
$$\le \sup_{n\in\N}k_n\frac{C^2+1}{C^2}E\left(\frac{\xi_{n1}^2}{1+\xi_{n1}^2}\right)<+\infty ; $$
By Theorem 2 in  [13] (p. 373),

$$\sup_{n\in\N}k_n E(\xi_{n1}I_{\{|\xi_{n1}|<C\}})<+\infty ; $$
$$\sup_{n\in\N}k_nD\left(\xi_{n1}I_{\{|\xi_{n1}|<C\}}\right)\le\sup_{n\in\N}k_nE\xi_{n1}^2I_{\{|\xi_{n1}|<C\}}=$$
$$=\sup_{n\in \NN}k_n(1+C^2)E\left(\frac{\xi_{n1}^2}{1+C^2}I_{\{|\xi_{n1}|<C\}}\right)\le \sup_{n\in\N}k_n(1+C^2)E\left(\frac{\xi_{n1}^2}{1+\xi_{n1}^2}\right) <+\infty.$$
The proof is completed.

 {\bf Theorem 5. } {\it Let (A) be valid and $X_n\stackrel{d}{\to}Y_\alpha$ as $n \to \infty$. Then it holds 
 \begin{equation}            \label{Qn}
 Q_n(\omega)\stackrel{w}{\to}\mu_\alpha
 \end{equation}
 as $n \to \infty$ for almost all $\omega \in \Omega$ in $D[0,1]$.}

 {\bf Proof. }
  Let $0<l<k$ and
 $$
 X_{lk}(t)=\left\{
 \begin{array}{lr}
b_kt& 0\leq t \leq \frac{s_l}{s_k}\\
X_{k}(t) - \frac{V(s_l)}{a_k}& \frac{s_l}{s_k}\leq t \leq 1\\
\end{array}
\right.
$$
Then $X_{lk} \mbox{ and } X_l$ are independent random processes.

It's clear, that $D[0,1]$ is a complete separable space in the metric  $$\rho_0(x,y)=\frac{\rho(x,y)}{1+\rho(x,y)}, \ \ x,y \in D[0,1]. $$
We will consider the metric
$$ \rho_1(x,y)=\min\{\sup_{0\le t\le 1}|x(t)-y(t)|,1\} .$$
Observe, that $$\rho_1(x,y)\ge \rho_0(x,y), \ \ x,\ y \ \in D[0,1] . $$

We will show that
\begin{equation}                    \label{12}
 E\rho_1(X_k,X_{kl})\le C\left(\frac{l}{k}\right)^{\frac{\beta}{2}}.
 \end{equation}
Since $V(t)$ is a random process with independent increments, 
$$ V(s_l)=\sum_{i=1}^{n_l}\xi_i, $$ where $\xi_i=V\left(\frac{s_li}{n_l}\right)-V\left(\frac{s_l(i-1)}{n_l}\right)$ are independent identically distributed  random variables.
 Let $n_k \in \N$ such that $\frac{n_l}{n_k}\le \frac{s_l}{s_k} $ and it exists $\lambda \in [0,1]$ such that 
 $$V(\lambda)=\sum_{i=1}^{n_k}\xi_i .$$
 
We have
$$E\min\left\{\frac{1}{a_k}\sum_{i=1}^{n_l}|\xi_iI_{\{|\xi_i|\ge a_k\}}|,1\right\}\le n_lE\min\{\frac{1}{a_k}|\xi_i|I_{\{|\xi_i|\ge a_k\}},1\}\le n_lP\{|\xi_1|\ge a_k\}=$$
$$=\frac{n_l}{n_k}n_kP\{|\xi_1|\ge a_k\}\le  C_1\left(\frac{l}{k}\right)^{\beta}\le C_1\left(\frac{l}{k}\right)^{\frac{\beta}{2}},$$ 
and
$$E\min\left\{\frac{1}{a_k}\left|\sum_{i=1}^{n_l}\xi_iI_{\{|\xi_i|<a_k\}}\right|,1\right\}=$$
$$E\min\{\frac{1}{a_k}|\sum_{i=1}^{n_l}\xi_iI_{\{|\xi_i|<a_k\}}-E\xi_iI_{\{|\xi_i|<a_k\}}+E\xi_iI_{\{|\xi_i|<a_k\}}|,1\}\le$$
$$\frac{s_l}{s_k}\sqrt{n_kD(\frac{1}{a_k}\xi_1I_{\{|\xi_1|<a_k\}})}+\frac{s_l}{s_k}n_kE(\frac{1}{a_k}\xi_1I_{\{|\xi_1|<a_k\}})\le (C_2+\sqrt{C_3})\left( \frac{l}{k}\right)^{\frac{\beta}{2}},$$
where $C_1, \ C_2, \ C_3 $ are defined in Lemma 1.

Then, using the moment inequality from [8], sect. 5, p. 231  we obtain
$$E\rho_1(X_k,X_{lk})=E\min\{\sup_{0\le t \le 1}|X_k-X_{lk}|,1\}= E\min\left\{\sup_{0\le t \le \frac{s_l}{s_k}}|\frac{V(ts_k)}{a_k}|,1\right\}\le$$ $$4E\min\left\{\left|\frac{V(s_l)}{a_k}\right|,1\right\}=4E\min\{\frac{1}{a_k}|\sum_{i=1}^{n_l}\xi_i|,1\} \le$$ $$4E\min\{\frac{1}{a_k}|\sum_{i=1}^{n_l}\xi_iI_{\{|\xi_i|<a_k\}}|,1\}+4E\min\{\frac{1}{a_k}|\sum_{i=1}^{n_l}\xi_iI_{\{|\xi_i|\ge a_k\}}|,1\}\le 4(C_1+C_2+\sqrt{C_3})\left(\frac{l}{k}\right)^{\frac{\beta}{2}}.$$

Then by Lemma 1 from [2], this  and Theorems 2-4 imply the convergence ($\ref{Qn}$). The proof is completed.
\bigskip

 \centerline{\bf  References.}

 \medskip

 [1] {\sc M. T. Lacey} and {\sc W. Philipp,} A note on almost sure central limit theorem,
{\it Statistics and Probability Letters}
 {\bf9}(2) (1990), 201--205.

 [2] {\sc A. Chuprunov} and {\sc I. Fazekas},
Almost sure versions of some analogues of the invariance principle,
{\it Publicationes Mathematicae, Debrecen}
{\bf 54}(3-4) (1999), 457--471.

[3] {\sc A. Chuprunov} and {\sc I. Fazekas}, Almost sure limit theorems for the Pearson statistic,
{\it Teor. Veroyatnost. i Primenen.} {\bf 48}(1), 162--169.

[4] {\sc A. Chuprunov} and {\sc I. Fazekas},
Integral analogues of almost sure limit theorems,
{\it Periodica Mathematica Hungarica}o
V.5 (1-2), 2005 pp. 61-78.

[5] {\sc B. V. Gnedenko} and {\sc A. N. Kolmogorov},
Limit Distributions for Sums of Independent Random Variables,
{\it  Addison-Wesley, Reading, Massachusetts},
1954.

[6] {\sc P. Billingsley}, Convergence of probability measures,
{\it John Wiley and Sons, New York}, 1968.

[7] {\sc D. Pollard}, Convergence of stochastic processes,
{\it Springer-Verlag, New York}, 1984.

[8] {\sc N. N. Vakhania}, {\sc V. I. Tarieladze} and {\sc S. A. Chobanian},
Probability distributions in Banach spaces,
{\it Nauka, Moscow}, 1985 (in Russian).

 [9] {\sc V. M. Kruglov} and {\sc V. Yu. Korolev},
Limit theorems for random sums,
{\it Moscow University Press, Moscow},
1990 (in Russian).

[10] {\sc I. I. Gihman, A. V Skorohod}, Theory of the random processes,
{\it Nauka, Moscow}, 1973 (in Russian).

[11] {\sc A. V. Skorohod}, Random processes with independent increments, {\it Probability theory and mathematical statistics, Nauka, Moscow}, 1964, (in Russian).

[12] {\sc V. M. Borodikhin}, Conditions for the compactness of subsets of measures on some metric spaces, {\it Teor. Veroyat. i Matem. Statist.} {\bf 3} (1970), 16-28, (in Russian). 

[13] {\sc A. N. Shiryaev}, Probability, {\it Nauka, Moscow}, 1997 (in Russian)
\end{document}